\documentclass[reqno,11pt]{amsart}
\usepackage{a4wide}

\usepackage{latexsym,amssymb,amsfonts}
\usepackage{amsmath}
\usepackage[pctex32]{graphics}
\usepackage{graphicx}

\newtheorem{thm}{Theorem}
\newcommand{\rr}{\mathbb{R}}
\renewcommand{\ss}{\mathbb{S}^{d-1}}
\newcommand{\ee}{\mathbb{E}}

\newcommand{\X}{\overline X}
\newcommand{\V}{\overline V}
\newcommand{\ds}{\displaystyle}

\begin{document}

\title{Mean-field limit for the stochastic Vicsek model}

\author{Fran\c cois Bolley}
\address{Universit\'e Paris-Dauphine, Ceremade, Place du Mar\'echal de Lattre de Tassigny, F-75775 Paris cedex 16, France}
\email{bolley@ceremade.dauphine.fr}

\author{Jos\'e A. Ca\~nizo}

\address{Departament de Mate\-m\`a\-tiques,
Universitat Aut\`onoma de Barcelona, E-08193 Bellaterra, Spain}
\email{canizo@mat.uab.es}

\author{Jos\'e A. Carrillo}

\address{ICREA (Instituci\'o Catalana de Recerca i
Estudis Avan\c cats) and Departament de Mate\-m\`a\-tiques,
Universitat Aut\`onoma de Barcelona, E-08193 Bellaterra, Spain}
\email{carrillo@mat.uab.es}

\date{February 1, 2011}

\maketitle


\begin{abstract}
  We consider the continuous version of the Vicsek model with noise,
  proposed as a model for collective behavior of individuals with a
  fixed speed.  We rigorously derive the kinetic mean-field partial
  differential equation satisfied when the number $N$ of particles
  tends to infinity, quantifying the convergence of the law of one
  particle to the solution of the PDE. For this we adapt a classical
  coupling argument to the present case in which both the particle
  system and the PDE are defined on a surface rather than on the whole
  space $\rr^d$. As part of the study we give existence and uniqueness
  results for both the particle system and the PDE.
\end{abstract}


\section*{Introduction}

The stochastic Vicsek model \cite{vicsek} arises in the study of
collective motion of animals and it is receiving lots of attention due
to the appearance of a phase transition \cite{CKJRF,chate}. A
continuum version and variants of this model have been proposed in the
recent works~\cite{DM08,DFL10}. Our objective is to rigorously derive
some continuum partial differential equations analysed in \cite{DM08}
from the stochastic Vicsek particle model.  This was carried out for a
family of collective behaviour models in~\cite{BCC10} following the
method of~\cite{Szn91}. The present models do not fall into this
analysis due to the evolution being defined on a surface as we explain
next.  In the models considered here, individuals are assumed to move
with a fixed cruising speed trying to average their orientations with
other individuals in the swarm in the presence of noise. This
orientation mechanism is modelled by locally averaging in space their
relative velocity to other individuals. More precisely, we are
interested in the behaviour of $N$ interacting $\rr^{2d}$-valued
processes $(X^{i}_t, V^{i}_t)_{t\geq 0}$ with $1 \leq i \leq N$ with
constant speed $\vert V^i_t \vert $, say unity. We define them as
solutions to the coupled Stratonovich stochastic differential
equations
\begin{equation}\label{visstra}
\begin{cases}
  d X_t^i = V_t^i \, dt, \\
\displaystyle{dV_t^i = \sqrt{2}\, P(V_t^i) \circ dB_t^i - P(V_t^i)
\left(\frac{1}{N} \sum_{j=1}^N K(X_t^i \! - \! X_t^j) (V_t^i - V_t^j) \right) dt.} \\
\end{cases}
\end{equation}
Here $P(v)$ is the projection operator on the tangent space at $v/|v|$ to
 the unit sphere $\ss$ in $\rr^d$, i.e.,
$$
P(v)= I - \frac{v \otimes v}{\vert v \vert^2}\, .
$$
This stochastic system is considered with independent and commonly
distributed initial data $(X^{i}_0, V^{i}_0) \in \rr^d \times \ss$
with $1 \leq i \leq N$. The $(B^i_t)_{t \geq 0}$ denote $N$
independent standard Brownian motions in $\rr^d$. The projection operator ensures
that $V^i_t$ keeps constant norm, equal to $1$. The second term in
the evolution of $V^i_t$ models the tendency of the particle $i$
to have the same orientation as the other particles, in a way
weighted by the interaction kernel $K$, as in the model proposed
by F.~Cucker and S.~Smale~\cite{CS2}. Let us observe that
$P(V_t^i) V_t^i =0$, so we can drop the corresponding term when writing
\eqref{visstra} to recover the usual formulations as in
\cite{DM08}.

We will work with stochastic processes defined on
$\rr^{2d}$ instead of $\rr^d\times\ss$. We will check later on
that solutions of \eqref{visstra} with initial data in
$\rr^d\times\ss$ remain there for all times. We have written
\eqref{visstra} in the Stratonovich sense, since the term involving noise
corresponds to  Brownian motion on the sphere
$\ss$ as in \cite[Section 1.4]{hsu} and \cite[Section V.31]{rw}.

By symmetry of the initial configuration and of the evolution, all
particles have the same distribution. Even though they are
initially independent, correlation builds up in time due to the
interaction term. Nevertheless, this interaction term is of order
$1/N,$ and thus, it seems reasonable that two of these interacting
particles (or a fixed number $k$ of them) become less and less
correlated as $N$ gets large (propagation of chaos).

Following~\cite{Szn91} we shall show that the $N$ interacting
processes $({X}^i_t,{V}^i_t)_{t \geq 0}$ respectively behave as
$N\to\infty$ like the auxiliary processes $(\X^i_t,\V^i_t)_{t \geq
0}$, solutions to
\begin{equation}
  \label{eq:nlSDE}
\begin{cases}
      d\X_t^i = \V_t^i \,dt,
      \\
     \displaystyle  d\V_t^i = \sqrt{2} \, P(\V_t^i)  \circ  dB_t^i - P(\V_t^i) (H * f_t) (\X_t^i,  \V_t^i) \, dt,
      \\
      (\X_0^i,\V_0^i)=(X_0^i,V_0^i), \quad f_t = \textrm{law} (\X_t^i, \V_t^i)
    \end{cases}
\end{equation}
in the Stratonovich sense. Here the Brownian motions $(B^i_t)_{t \geq 0}$ are those governing the
evolution of the $(X^i_t,V^i_t)_{t \geq 0}$ and
$$
(H * f) (x,v) = \int_{\rr^{2d}} K(x-x') \, (v-v') \, f(x',v') \,
dx' \, dv', \qquad x, v \in \rr^d \,.
$$
Note that~\eqref{eq:nlSDE} consists of $N$ equations which can be
solved independently of each other. Each of them involves the
condition that $f_t$ is the distribution of $(\X^i_t,\V^i_t)$,
thus making it nonlinear. The processes ${(\X^i_t,\V^i_t)_{t \geq
0}}$ with $i\geq 1$ are independent since the initial conditions
and driving Brownian motions are independent.

We will show that these processes defined on $\rr^{2d}$ are
identically distributed, take values in $\rr^d \times \ss$ if
initially so, and their common law $f_t$ at time $t$, as a measure
on $\rr^d \times \ss$, evolves according to
\begin{equation}\label{eq:pde}
  \partial_t f_t + \omega  \cdot \nabla_x f_t  = \Delta_{\omega} f_t + \nabla_{\omega} \cdot \big( f_t (I - \omega \otimes \omega) (H*f_t) \big), \quad t>0, x \in \rr^d, \omega \in \ss.
\end{equation}
Now the convolution $H*f$ is over $\rr^{d} \times \ss$:
$$
(H*f) (x,\omega) = \int_{\rr^{d} \times \ss} K(x-x') \, (\omega -
\omega') \, f(x',\omega') \, dx' \, d\omega', \qquad x \in \rr^d,
\omega \in \ss \, .
$$
Moreover, $\nabla_x$ stands for the gradient with respect to the
position variable $x \in \rr^d$ whereas $\nabla_{\omega}$,
$\nabla_{\omega} \cdot$ and $\Delta_{\omega}$ respectively stand
for the gradient, divergence and Laplace-Beltrami operators with
respect to the velocity variable $\omega \in \ss$.

This equation is proposed in \cite{DFL10} as a continuous version of
the original Vicsek model \cite{vicsek}, and one of our purposes is to
make this derivation rigorous. The asymptotic behavior and the
appearance of a phase transition in the space-homogeneous version
of~\eqref{eq:pde} (i.e., without the space variable) has been recently
studied in \cite{FL11}. It is also known as the Doi-Onsager equation,
introduced by Doi in \cite{Doi} as a model for the non-equilibrium
Statistical Mechanics of a suspension of polymers in which their
spatial orientation (given by the parameter $\omega \in \ss$) is taken
into account.

\medskip

The main result of this paper can be summarized as:

\begin{thm}\label{thm:exuniq}
  Let $f_0$ be a probability measure on $\rr^d \times \ss$ with finite
  second moment in $x \in \rr^d$ and let $(X^i_0, V^i_0)$ for $1 \leq
  i \leq N$ be $N$ independent variables with law $f_0$. Let also $K$
  be a Lipschitz and bounded map on $\rr^d$. Then,
  \begin{enumerate}
  \item[i)] There exists a pathwise unique global solution to the SDE
    system~\eqref{visstra} with initial data $(X^i_0, V^i_0)$ for $1
    \leq i \leq N$; moreover, the solution is such that all $V^i_t$
    have norm $1$.
  \item[ii)] There exists a pathwise unique global solution to the
    nonlinear SDE~\eqref{eq:nlSDE} with initial datum $(X^i_0,
    V^i_0)$; moreover, the solution is such that $\V^i_t$ has norm
    $1$.
  \item[iii)] There exists a unique global weak solution to the
    nonlinear PDE \eqref{eq:pde} with initial datum $f_0$. Moreover,
    it is the law of the solution to \eqref{eq:nlSDE}.
  \end{enumerate}
\end{thm}
Solutions to general SDE's can be built in submanifolds of $\rr^d$ by
means of the Brownian motion of the ambient space as in \cite[Theorem
V.34.86]{rw} for instance; then one can interpret the generator in
terms of the corresponding Laplace-Beltrami operator. For example, the
Brownian motion on a submanifold $\Sigma$ of $\rr^d$ is the solution
to the SDE
$$
dW_t=P_\Sigma(W_t)\circ dB_t
$$
on $\rr^d$ and with $P_\Sigma (w)$ being the orthogonal projection of
$\rr^d$ onto the tangent space at $w$ to $\Sigma$. Here, we give the
full construction and derivation of the evolution of the law as it can
be done explicitly in the case of the sphere $\ss$.  Let us also
emphasize that we have only partial diffusion since it is a kinetic
model.

We observe that existence of $L^2$ and classical solutions for the
space-homogeneous version of \eqref{eq:pde} has also been considered in
\cite{FL11}.

\medskip

As a direct consequence of the classical Sznitman's theory, we
get the following mean-field limit result:

\begin{thm}\label{thm:smf}
With the assumptions of Theorem~\ref{thm:exuniq} and for the
respective solutions
 $(X^i_t, V^i_t)_{t \geq 0}$ and $(\X^i_t, \V^i_t)_{t \geq 0}$ of~\eqref{visstra} and~\eqref{eq:nlSDE}, for all $T >0$ there exists a constant $C$ such that
$$
 \ee \Big[ \vert X^i_t - \X^i_t \vert^2 + \vert V^i_t - \V^i_t \vert^2 \Big] \leq \frac{C}{N}
$$
for all $0\leq t \leq T$, $N \geq 1$ and $1 \leq i \leq N$.
\end{thm}

This estimate classically ensures quantitative estimates on (see
\cite{Szn91,BCC10} for details)
\begin{enumerate}
  \item[i)]
the convergence in $N$ of the law at time $t$ of any (by symmetry)
of the processes $(X^i_t, V^i_t)$ towards $f_t$,
  \item[ii)]
the propagation of chaos for the particle system through the
convergence of the law at time $t$ of any $k$ particles towards
the tensor product $f_t^{\otimes k}$ (for $k$ fixed or $k=o(N)$),
  \item[iii)]
the convergence of the empirical measure at time $t$ of the
particle system towards $f_t$.
\end{enumerate}
Of course, the same techniques lead to a corresponding mean-field
limit result for the space-homogeneous particle system instead of
\eqref{visstra}, obtaining the corresponding space-homogeneous PDE.

\section*{Proofs}

Using the standard It\^o-Stratonovich calculus, see~\cite[p.
99]{gardiner} for instance, equations~\eqref{visstra} and~\eqref{eq:nlSDE}
are respectively equivalent to the It\^o stochastic
differential equations
\begin{equation}
\label{eq:sdesysITO}
\begin{cases}
  d X_t^i = V_t^i dt, \\
\displaystyle{dV_t^i = \sqrt{2} \, P(V_t^i) dB_t^i - P(V_t^i) \left( \frac{1}{N} \sum_{j=1}^N
K(X_t^i \! - \! X_t^j) (V_t^i - V_t^j) \right) dt - (d-1) \frac{V^i_t}{\vert V^i_t\vert^2} \, dt} \,.\\
\end{cases}
\end{equation}
and
\begin{equation}
  \label{eq:nlSDEITO}
\begin{cases}
      d\X_t^i = \V_t^i \,dt,
      \\
     \displaystyle  d\V_t^i = \sqrt{2} \, P(\V_t^i) dB_t^i - P(\V_t^i) (H * f_t) (\X_t^i,  \V_t^i) \, dt - (d-1) \frac{\V^i_t}{\vert \V^i_t\vert^2} \, dt,
      \\
      (\X_0^i,\V_0^i)=(X_0^i,V_0^i), \quad f_t = \textrm{law} (\X_t^i, \V_t^i)
    \end{cases}
\end{equation}
which we now consider.

\medskip

We start with the proof of Theorem~\ref{thm:exuniq}. We use a regularization of the diffusion and
drift coefficients. We let $\sigma_1$ be a $d \times d$ matrix valued
map on $\rr^d$ with bounded derivatives of all orders such that
$\sigma_1 (v) = P(v)$ for all $v$ with $\vert v \vert \geq 1/2$, and
$\sigma_2$ and $\sigma_3$ be maps on $\rr^d$, again with bounded
derivatives of all orders, such that $\sigma_2(v) = v / \vert v
\vert^2$ if $\vert v \vert \geq 1/2$ and $\sigma_3 (v) = v$ if $\vert
v \vert \leq 2$.

\subsection*{Existence and uniqueness for the particle system~\eqref{eq:sdesysITO}}\label{subsec:partsys}

Given such $\sigma_1, \sigma_2$, the system of equations
\begin{equation}
\label{eq:sdesyssigma}
\begin{cases}
  d X_t^i = V_t^i dt, \\
\displaystyle{dV_t^i = \sqrt{2} \, \sigma_1(V_t^i) \, dB_t^i -
\sigma_1(V_t^i) \!\left( \frac{1}{N} \sum_{j=1}^N K(X_t^i \!-\!
X_t^j) (V_t^i - V_t^j) \right) \! \! dt - (d-1) \sigma_2(V^i_t) dt}
\end{cases}
\end{equation}
starting from $(X^i_0, V^i_0) \in \rr^d \times \ss$ for $1 \leq i
\leq N$ has locally Lipschitz coefficients. Moreover, by the It\^o
formula and as long as $\vert V^i_t \vert \geq 1/2$,
\begin{eqnarray*}
d \vert V^i \vert^2
&=&
2 \sqrt{2} \, V^i \cdot P(V^i)  dB^i  - 2 \, V^i \cdot P(V^i) \left( \frac{1}{N} \sum_{j=1}^N K(X^i - X^j) (V^i - V^j) \right) dt \\
&& - 2 \, (d-1)  dt + 2 \sum_{k, l=1}^d  \delta_{kl} \; d \left< B^i_k - \sum_{p=1}^d \frac{V^i_k V^i_p}{\vert V^i \vert^2} B^i_p \, , \, B^i_l - \sum_{q=1}^d \frac{V^i_l V^i_q}{\vert V^i \vert^2} B^i_q \right>\\
&=& -2 \, (d-1) dt+ 2 \sum_{k=1}^d \left[ 1 - 2 \, \frac{(V^i
_k)^2}{\vert V^i \vert^2} + \sum_{p=1}^d \frac{(V^i_p)^2
(V^i_k)^2}{\vert V^i \vert^4} \right] dt = 0.
\end{eqnarray*}
Here we dropped the time dependence, wrote $y = (y_1, \dots y_d)
\in \rr^d$ and used the fact that $\ds V^i \cdot P(V^i) y =0$ for
all vectors $y\in \rr^d$. Hence $\vert V^i_t \vert =1$ up to
explosion time. Since moreover $dX^i_t = V^i_t \, dt$, this
ensures that the explosion time is infinite, hence global
existence and pathwise uniqueness for~\eqref{eq:sdesyssigma}.

Now the solution to~\eqref{eq:sdesyssigma} for given $\sigma_1,
\sigma_2$ is a solution to~\eqref{eq:sdesysITO} since all velocities
have norm~$1$, which provides global existence of solutions
to~\eqref{eq:sdesysITO}. If now we consider two solutions
to~\eqref{eq:sdesysITO} for the same initial data and Brownian
motions, then they have velocities equal to $1$, so that are
solutions to~\eqref{eq:sdesysITO} for any $\sigma_1, \sigma_2$, for
which pathwise uniqueness holds: hence they are equal. This proves
the first part in Theorem~\ref{thm:exuniq}.

\subsection*{Existence and uniqueness for the artificial processes~\eqref{eq:nlSDEITO}}
Let $\sigma_1, \sigma_2, \sigma_3$ be any maps as above and let
$$
H_{\sigma_3} [f] (x) = \int_{\rr^{2d}} K(x-y) \, \sigma_3(v-w) \,
f(y,w) \, dy \, dw.
$$
Then, given a distribution $f_0$ on $\rr^d \times \ss$ with finite
second moment in $x \in \rr^d$ and $(\X_0,\V_0)$ with law $f_0$,
the nonlinear equation
\begin{equation}
  \label{eq:sdenl1}
  \begin{cases}
    d\X_t = \V_t \,dt,
    \\
    \displaystyle  d\V_t= \sqrt{2}\, \sigma_1(\V_t) dB_t
    - \sigma_1(\V_t)  (H_{\sigma_3} * f_t) (\X_t, \V_t) dt - (d-1) \sigma_2(\V_t) dt,
    \\
    f_t = \textrm{law} (\X_t, \V_t)
  \end{cases}
\end{equation}
has bounded and Lipschitz coefficients on $\rr^{2d}$, so admits a
pathwise unique global solution according to~\cite[Theorem
1.1]{Szn91}. Moreover, as long as $ \vert \V_t  \vert \geq 1/2$,
then we can repeat the argument above  to
prove that $d \vert \V_t \vert^2 = 0$, so that $\vert \V_t \vert
=1$ for all time. In particular the obtained solution
$(\X_t,\V_t)_{t \geq 0}$ is a global solution to the genuine
nonlinear equation~\eqref{eq:nlSDEITO}. Pathwise uniqueness of
solutions to~\eqref{eq:nlSDEITO} can be obtained as
for~\eqref{eq:sdesysITO}.

\subsection*{Existence and uniqueness for the PDE~\eqref{eq:pde}}
Let $f_0$ be a distribution on $\rr^d \times \ss$ with finite
second moment in $x \in \rr^d$, $(\X_0,\V_0)$ with law $f_0$, and
let $(\X_t,\V_t)_{t \geq 0}$ be the solution to~\eqref{eq:nlSDEITO}
with initial datum $(\X_0,\V_0)$. Its law $f_t$, as a measure on
$\rr^{2d}$, satisfies
$$
\frac{d}{dt} \int_{\rr^{2d}} \varphi \, df_t = \int_{\rr^{2d}}
\big( v \cdot \nabla_x \varphi  + \textrm{Hess}_v \varphi : (I - v
\otimes v) +  \nabla_v \varphi \cdot (I - v \otimes v) (H*f_t) -
(d-1) v \cdot \nabla_v \varphi   \big) df_t
$$
for all smooth $\varphi$ on $\rr^{2d}$ by the It\^o formula; here
$\nabla_v$ and $\Delta_v$ are respectively the gradient and
Laplace operators with respect to $v \in \rr^d$, and
$\textrm{Hess}_v \varphi : M$ is the term by term product of the
Hessian with respect to $v$ matrix of $\varphi$ with a matrix $M$.

We have observed that $\vert \V_t \vert =1$ {\it a.s.}, so $f_t$
is concentrated on $\rr^d \times \ss$. We now define the
restriction $F_t$ of $f_t$ on $\rr^d \times \ss$ by
$$
\int_{\rr^d \times \ss} \Phi \, dF_t = \int_{\rr^{2d}}  \varphi \, df_t
$$
for all continuous maps $\Phi$ on $\rr^d \times \ss$, where
$\varphi$ is any continuous and bounded map on $\rr^{2d}$ equal to
$\Phi$ on $\rr^d \times \ss$. Let now $\Phi$ be a $\mathcal
C_c^{\infty}$ map on $\rr^d \times \ss$ and $\varphi$ be a $\mathcal
C^{\infty}_c$ map on $\rr^{2d}$ such that $\varphi(x, v) = \Phi(x,
v / \vert v \vert)$ for all $1/2 \leq \vert v \vert \leq 2$. Then
$\varphi$ is $0$-homogeneous in $v$ in the annulus $1/2 \leq \vert
v \vert \leq 2$, so that $v \cdot \nabla_v \varphi = 0$ for all
$(x,v)$ in the support of $f_t$. In particular
$$
\frac{d}{dt} \int_{\rr^{d} \times \ss} \Phi \, dF_t = \frac{d}{dt}
\int_{\rr^{2d}} \varphi \, df_t = \int_{\rr^{2d}} \Big( v \cdot
\nabla_x \varphi  + \Delta_v \varphi +  \nabla_v \varphi \cdot (I
- v \otimes v) (H*f_t)   \Big) df_t.
$$
Then the maps $v \cdot \nabla_x \Phi$ and $v \cdot \nabla_x
\varphi$ are equal on $\rr^d \times \ss$ since $\Phi$ and
$\varphi$ have the same $x$-dependence. Moreover, $\nabla_\omega
\Phi = \nabla_v \varphi$ and $\Delta_{\omega} \Phi = \Delta_v
\varphi$ for $(x, \omega) \in \rr^d \times \ss$. This last point
can be checked by direct computations. Hence
$$
\frac{d}{dt} \int_{\rr^{d} \times \ss} \Phi \, dF_t =
\int_{\rr^{d} \times \ss} \left( \omega \cdot \nabla_x \Phi  +
\Delta_{\omega} \Phi + \nabla_{\omega} \Phi \cdot (I - \omega
\otimes \omega) (H*F_t)  \right)  dF_t.
$$
This ensures that $F_t$ is a weak solution to~\eqref{eq:pde}.

\medskip

We now turn to uniqueness of solutions to~\eqref{eq:pde}. For that
purpose we let $f^1$ and $f^2$ be two solutions with the same initial
datum $f_0$, and at each time $t$ we view them as measures on
$\rr^{2d}$ concentrated on the surface $\rr^d \times \ss$. We let
$(\X^1_t, \V^1_t)_{t \geq 0}$ and $(\X^2_t, \V^2_t)_{t \geq 0}$ be the
solutions to~\eqref{eq:sdenl1} with drift given by $H_{\sigma_3} *
f_t^1$ and $H_{\sigma_3} * f_t ^2$ respectively, and common initial
datum $(\X_0, \V_0)$ with law $f_0$. Then their respective laws
$g^1_t$ and $g^2_t$, as measures on $\rr^{2d}$, are solutions to the
linear PDE
\begin{equation*}
  \partial_t g^i_t + v  \cdot \nabla_x g^i_t
  = \sum_{k,l=1}^d \frac{\partial^2}{\partial v_k \partial v_l} \left(
    (\sigma_1 \sigma_1^*)_{kl} \, g^i_t
  \right)
  + \nabla_v \cdot \left[
    g^i_t \left(
      \sigma_1\, (H_{\sigma_3} * f^i_t) + (d-1) \sigma_2
    \right)
  \right].
\end{equation*}
Since $f^i_t$ is also a measure solution to this linear PDE on
$\rr^{2d}$ with bounded and regular coefficients, for which uniqueness
classically holds, it follows that $g^i_t = f^i_t$
($i=1,2$). Consequently, the $(\X^i_t, \V^i_t)_{t \geq 0}$ are
solutions to the nonlinear SDE~\eqref{eq:sdenl1}, for which we have
already proved uniqueness. Hence $(\X^1_t, \V^1_t)_{t \geq 0}$ and
$(\X^2_t, \V^2_t)_{t \geq 0}$ are equal, and in particular $f^1_t (=
g^1_t) = (g^2_t=) f^2_t$.

\subsection*{Proof of Theorem~\ref{thm:smf}}

Since $\vert V^i _t \vert = \vert \V^i_t \vert =1$  for all $i$
and $t$, the processes $(X^i_t, V^i_t)_{t \geq 0}$ and $(\X^i_t,
\V^i_t)_{t \geq 0}$ are solutions of the corresponding equations
with bounded and Lipschitz diffusion and drifts coefficients as
in~\eqref{eq:sdenl1}. Hence we may apply the estimates
in~\cite[Theorem 1.4]{Szn91} to obtain Theorem~\ref{thm:smf}.

\bigskip

\footnotesize \noindent\textit{Acknowledgments.} The authors wish
to thank Pierre Degond and Vlad Panferov for drawing their
attention to this question and for discussion on it. This note was
written while the authors were visiting the Isaac Newton
Institute, Cambridge; it is a pleasure for them to thank this
institution for its kind hospitality.

The last two authors acknowledge support from the project
MTM2008-06349-C03-03 DGI-MCI (Spain), the 2009-SGR-345 from
AGAUR-Generalitat de Catalunya and the French-Spanish acciones
integradas program FR2009-0019. All authors were partially
supported by the ANR-08-BLAN-0333-01 Projet CBDif-Fr.

\end{document}